\documentclass[11pt]{article}
\usepackage{amssymb}
\usepackage{amsthm,undertilde}
\usepackage{amsmath}
\usepackage{amsbsy}
\usepackage{epsfig}
\usepackage{color}
\usepackage{enumerate}
\usepackage{vmargin}
\setpapersize{USletter}
\setmarginsrb{1in}{1in}{1in}{1in}{0pt}{0mm}{0pt}{36pt}

\begin{document}

\begin{center}

{\bf \Large Intrinsic Random Functions on the sphere}

Chunfeng Huang\footnote{Department of Statistics, Indiana University, Bloomington, IN. Email: huang48{@}indiana.edu.}, 
Haimeng Zhang\footnote{Department of Mathematics and Statistics, University of North Carolina at Greensboro, Greensboro, NC. Email:haimengzhanguncg{@}gmail.com}, 
Scott M. Robeson\footnote{Department of Geography, Indiana University, Bloomington, IN. Email:srobeson{@}indiana.edu}, 
Jacob Shields\footnote{Department of Statistics, Indiana University, Bloomington, IN. Email:shields{@}imail.iu.edu}

\end{center}

{\bf Abstract.} Spatial stochastic processes that are modeled over the entire Earth's surface require statistical approaches that directly consider the spherical domain. Here, we extend the notion of intrinsic random functions (IRF) to model non-stationary processes on the sphere and show that low-frequency truncation plays an essential role. Then, the universal kriging formula on the sphere is derived. We show that all of these developments can be presented through the theory of reproducing kernel Hilbert space. In addition, the link between universal kriging and splines is carefully investigated, whereby we show that thin-plate splines are non-applicable for surface fitting on the sphere.

{\bf Keywords.} Kriging, Reproducing kernel Hilbert space, Splines, Stationary

\section{Introduction}

Global-scale phenomena, of which there are a multitude of important applications ranging from climate science to epidemiology, can be viewed as random processes on the sphere. A common assumption when modeling such processes is second-order stationarity (or stationarity for short in this paper), where the mean is constant and the covariance function at two locations is assumed to depend only on their distance (or, in the case of circular or spherical phenomena, angular separation; see \cite{Obukhov1947, Huangetal2011, Gneiting2013}). This assumption is  difficult to evaluate and often deemed unrealistic in practice. Several approaches have been proposed to relax this assumption such as axial symmetry \cite{Jones1963, Stein2007, Huangetal2012} and kernel convolution \cite{ZhuandWu2010, Heatonetal2014}. Based on generalized random functions in \cite{GelfandandVilenkin1964}, Matheron \cite{Matheron1973} introduced a flexible family of non-stationary processes, named intrinsic random functions (IRFs), where the process is assumed to have lower monomials as its mean and the transformed process becomes stationary. The influence of this approach is profound, as it provides a basis for kriging \cite{Cressie1993, Stein1999, ChilesandDelfiner2012} in practice. Most developments and uses of IRFs have been in Euclidean spaces, where the differencing or differential operation is used for the transformation. For example, IRFs on the real line, also known as processes with random stationary increments, were discussed in \cite{Yaglom1955} and \cite{ZhangandHuang2014}.  IRFs on other spaces are less explored. In Huang et al. \cite{Huangetal2016}, IRFs on the circle are developed, where we showed that low-frequency truncation replaces differencing to achieve stationarity. 

In this paper, we first define IRFs on the sphere. While lower monomials are used in Euclidean spaces, the counterparts on the sphere are the lower order spherical harmonics. We formally show that an IRF on the sphere is characterized by its frequency-truncated process (Theorem 1 in Section 2.1). Then, we derive the universal kriging formula for IRFs on the sphere, where we show that the coefficients for lower spherical harmonics do not need to be estimated. In addition, we demonstrate that IRFs can be viewed through the theory of reproducing kernel Hilbert space (RKHS, \cite{Aronszajn1950, Wahba1990a, Taijeronetal1994, Levesleyetal1999, Huangetal2016}). Based on this, kriging on the sphere is shown to be the same as the smoothing formula in RKHS. We formally establish the equivalence between splines and kriging. Finally, we carefully investigate splines on the sphere and find that the popularly used thin-plate spline approach is not applicable for surface fitting on the sphere.  \\

\section{Main results}

Let $Z(x)$ be a real-valued random process on a unit sphere $S^2$, where $x=(\psi,\zeta)$ with the longitude $\psi \in [0, 2\pi)$ and the latitude $\zeta \in [0, \pi]$. Assume $Z(x)$ is continuous in quadratic mean and its covariance function $\mbox{cov} (Z(x), Z(y))$ is strictly positive definite. One can expand $Z(x)$ in spherical harmonics which are convergent in quadratic mean \cite{Yaglom1961, Jones1963, Roy1969, Roy1973}:
\[
Z(x) = \sum_{l=0}^\infty \sum_{m=-l}^l Z_{l,m} Y_l^m(x),
\]
where the real-valued spherical harmonic functions are
\[
\left\{ \begin{array}{ll} Y_l^m(x) = \sqrt{ \frac{2l+1}{2\pi} \frac{(l-m)!}{(l+m)!}} P_l^m(\cos \zeta) \cos (m \psi), & m=1,\ldots, l,\\
Y_l^0(x) = \sqrt{ \frac{2l+1}{4\pi}} P_l^0(\cos \zeta), & \\
Y_l^{-m}(x) = \sqrt{ \frac{2l+1}{2\pi} \frac{ (l-m)!}{(l+m)!}} P_l^m(\cos \zeta) \sin (m \psi), & m=1,\ldots, l, \end{array} \right.
\]
and $P_l^m(\cdot)$ are the associated Legendre polynomials and $P_l^0(\cdot) \equiv P_l(\cdot)$ are the Legendre polynomials. The coefficients $Z_{l,m}$ are random variables: 
\[
Z_{l,m} = \int_0^{2\pi} \int_0^{\pi} Z(\psi, \zeta) Y_l^m(\psi, \zeta) \sin \zeta d\zeta d\psi.
\]
We denote the group of rotations of the sphere by $\mathcal{G}$. When the process $Z(x)$ is assumed to be stationary, that is, the mean is constant $\mbox{E} (Z(x)) = \mu$ and, for any two points $x, y \in S^2$,
\[
\mbox{cov} (Z(x), Z(y)) = \mbox{cov} (Z(gx), Z(gy)),
\]
for all $g \in \mathcal{G}$. The stationarity implies that the coefficients $Z_{l,m}$ are uncorrelated \cite{Obukhov1947, Yaglom1961, Roy1969}:
\[
\mbox{cov} (Z_{l,m}, Z_{l',m'}) = a_l I(l, l') I(m,m'),
\]
where $a_l \ge 0, \sum_l (2l+1) a_l < \infty$ and $I(l,l')$ is the indicator function taking value $1$ if $l=l'$, and zero otherwise. By the addition theorem of Legendre polynomials, the covariance function
\begin{equation} \label{eq:cov}
\mbox{cov} (Z(x), Z(y)) = \sum_{l=0}^\infty a_l \sum_{m=-l}^l Y_l^m(x) Y_l^m(y) = \sum_{l=0}^\infty \frac{2l+1}{4\pi} a_l P_l(\cos (d(x,y))),  
\end{equation}
where $d(x,y)$ is the spherical distance (angle) between $x=(\psi_x,\zeta_x)$ and $y=(\psi_y, \zeta_y)$:
\begin{equation} \label{eq:dxy}
d(x,y) = \cos^{-1} \left\{ \sin \zeta_x \sin \zeta_y + \cos \zeta_x \cos \zeta_y \cos (\psi_x-\psi_y) \right\}.
\end{equation}
The spectral representation of a positive definite function in the form of (\ref{eq:cov}) is also given in \cite{Schoenberg1942}. \\

\subsection{Intrinsic random functions}

Let $C(S^2)$ be the space of real-valued continuous functions on the sphere. By the Riesz representation theorem (\cite{Royden1988}, Chapter 13), the dual of $C(S^2)$ is the set of all finite regular signed Borel measures on $S^2$, denoted by $\Lambda$. Hence, for $f \in C(S^2)$ and $\lambda \in \Lambda$, we define
\[
f(\lambda) = \int_{S^2} f(x) \lambda(dx),
\]
and
\[
f(g \lambda) = \int_{S^2} f(gx) \lambda(dx).
\]
for a rotated measure $g\lambda$, where $ g \in \mathcal{G}$. 

Given an integer $\kappa > 0$, a measure $\lambda \in \Lambda$ is called an allowable measure of order $\kappa$  if it annihilates the spherical harmonics of order less than $\kappa$, that is, 
\[
\int_{S^2} Y_l^m(x) \lambda(dx) = 0, \quad 0 \le l < \kappa, \quad |m| \le l.
\]
Let $\Lambda_{\kappa}$ be the class of all allowable measures of order $\kappa$. It is obvious that $\Lambda_{\kappa+1} \subset \Lambda_{\kappa} \subset \Lambda$. Next, we establish the invariance of such allowable measures. \\

\noindent
{\bf Proposition 1.} For $\lambda \in \Lambda_{\kappa}$, the rotated measure $g\lambda$ is invariant, that is, $g \lambda \in \Lambda_{\kappa}$.
\begin{proof}
For any $g \in \mathcal{G}$, we have (\cite{Roy1969} and \cite{MarinucciandPeccati2011}, Ch. 3), 
\[
Y_l^m(gx) = \sum_{m'=-l}^l S(l,m,m') Y_l^{m'} (x), \quad x \in S^2, 
\]
where $S(l) = \{ S(l,m,m') \}_{m,m'=-l,\ldots, l}$ is a $(2l+1) \times (2l+1)$ orthogonal matrix that is independent of $x$. Let $\lambda \in \Lambda_{\kappa}$, then $Y_l^m(\lambda) = 0$ for $l < \kappa$ and $|m| \le l$. This leads to 
\[
Y_l^m(g\lambda) = \int_{S^2} Y_l^m(gx) \lambda(dx) = \sum_{m'=-l}^l S(l,m,m') \int_{S^2} Y_l^{m'}(x) \lambda(dx) = 0.
\]
\end{proof}

Proposition 1 shows that an allowable measure is invariant with respect to rotations of the sphere, hence the set $\Lambda_{\kappa}$ is closed for rotations. In short, 
\[
Y_l^m(\lambda)=0 \mbox{ and } Y_l^m(g\lambda)=0, \quad \mbox{for } \lambda \in \Lambda_{\kappa}, l < \kappa, |m| \le l, \mbox{ and } g \in \mathcal{G}.
\]
Parallel results on the circle are obtained in \cite{Huangetal2016}. Such invariance is essential in developing IRFs, see \cite{Matheron1973, Matheron1979}.   \\

Now, we extend the definitions of $f(\lambda)$ and $f(g\lambda)$ to the continuous random process $Z(x)$. If $Z(x)$ is with the probability measure space $(\Omega, \mathcal{A}, P)$, we denote
\[
Z(\lambda) = \int_{S^2} Z(x) \lambda(dx),
\]
that is, $Z(\lambda)$ maps $\Lambda$ to $L^2(\Omega, \mathcal{A}, P)$, and
\[
Z(g \lambda) = \int_{S^2} Z(gx) \lambda(dx).
\]
One can define the inner product
\[
\langle Z(\lambda_1), Z(\lambda_2) \rangle = \int_{S^2} \int_{S^2} \mbox{cov}(Z(x), Z(y)) \lambda_1(dx) \lambda_2(dy).
\]
Then, $\lambda \to Z(\lambda)$ is a continuous mapping, and if $\lambda_1=\lambda_2$, $Z(\lambda_1)=Z(\lambda_2)$ in the sense of norm induced by this inner product. \\

\noindent
{\bf Definition.} For an integer $\kappa$, a continuous random process $Z(x)$ on the sphere is called an Intrinsic Random Function of order $\kappa$ (IRF$\kappa$), if for any $\lambda \in \Lambda_{\kappa}$, the process is invariant with respect to any $g \in \mathcal{G}$. That is, $\mbox{E} (Z(\lambda)) = \mbox{E} (Z(g\lambda))$, and 
\[
\mbox{cov} (Z(\lambda_1), Z(\lambda_2)) = \mbox{cov} (Z(g \lambda_1), Z(g \lambda_2)), \quad \lambda_1, \lambda_2 \in \Lambda_{\kappa}.
\]

The following theorem characterizes the IRF on the sphere and reveals its connection to the stationary process through low-frequency truncation.\\

\noindent
{\bf Theorem 1.} For an integer $\kappa \ge 1$, a continuous random process $Z(x)$ on the sphere is an IRF$\kappa$ if and only if its frequency-truncated process $Z_{\kappa}(x)$ is stationary, where
\[
Z_{\kappa}(x) = \sum_{l=\kappa}^\infty \sum_{m=-l}^l Z_{l,m} Y_l^m(x).
\]

\begin{proof}
``$\Leftarrow$", for an integer $\kappa \ge 1$, assume the frequency-truncated process $Z_{\kappa}(x)$ is stationary. Note that 
\[
Z(x) = \sum_{l < \kappa} \sum_{m=-l}^l Z_{l,m} Y_l^m(x) + Z_{\kappa}(x),
\]
and $Z(\lambda) = \int_{S^2} Z(x) \lambda(dx)$ for $\lambda \in \Lambda_{\kappa}$. We have
\begin{eqnarray*}
&\,& Z(\lambda) = \int_{S^2} \left( \sum_{l < \kappa} \sum_{m=-l}^l Z_{l,m} Y_l^m(x)  + Z_{\kappa} (x) \right) \lambda (dx) \\
&\,& = \int_{S^2} \left( \sum_{l < \kappa} \sum_{m=-l}^l Z_{l,m} Y_l^m(x) \right) \lambda(d x) + \int_{S^2} Z_{\kappa}(x) \lambda(dx).
\end{eqnarray*}
The first term is with finite summation of $Y_l^m(\lambda)$ and it vanishes since $\lambda \in \Lambda_{\kappa}$. Therefore,
\[
Z(\lambda) = \int_{S^2} Z_{\kappa}(x) \lambda(dx) = Z_{\kappa} (\lambda), \quad \lambda \in \Lambda_{\kappa}.
\]
The expectation of $Z_{\kappa}(\lambda)$ is
\[
\mbox{E} (Z_{\kappa}(\lambda)) = \mbox{E} \int_{S^2} Z_{\kappa}(x) \lambda (dx).
\]
Since $\lambda$ is a finite measure and $Z_{\kappa}(x)$ is continuous in quadratic mean on the sphere, by Fubini's Theorem, we exchange the expectation and integral and obtain 
\[
\mbox{E} (Z_{\kappa}(\lambda)) = \int_{S^2} \mbox{E} (Z_{\kappa}(x)) \lambda(dx).
\]
Similarly, for $\lambda \in \Lambda_{\kappa}$ and $g \in \mathcal{G}$,
\[
Z(g\lambda) = Z_{\kappa}(g \lambda), \quad \mbox{and} \quad \mbox{E}(Z_{\kappa}(g\lambda)) = \int_{S^2} \mbox{E} (Z_{\kappa}(gx)) \lambda(dx).
\]
Then, following that $Z_{\kappa}(x)$ is stationary with $\mbox{E}(Z_{\kappa}(x)) = \mbox{E} (Z_{\kappa}(gx))$, we have  $\mbox{E} (Z_{\kappa}(g\lambda)) = \mbox{E} (Z_{\kappa}(\lambda))$ and arrive at
\[
\mbox{E} (Z(\lambda)) = \mbox{E} (Z(g\lambda)), \quad \lambda \in \Lambda_{\kappa}. 
\]

For the covariance part, when $\lambda_1, \lambda_2 \in \Lambda_{\kappa}$, we have $\mbox{cov} (Z(\lambda_1), Z(\lambda_2)) = \mbox{cov} (Z_{\kappa}(\lambda_1), Z_{\kappa}(\lambda_2))$ and $\mbox{cov} (Z(g\lambda_1), Z(g \lambda_2))=\mbox{cov} (Z_{\kappa}(g \lambda_1), Z_{\kappa}(g\lambda_2))$.  Similarly, $Z_{\kappa}(x)$ is stationary, which implies $\mbox{cov} (Z_{\kappa}(\lambda_1), Z_{\kappa}(\lambda_2)) = \mbox{cov} (Z_{\kappa}(g \lambda_1), Z_{\kappa}(g \lambda_2))$. This leads to
\[
\mbox{cov} (Z(\lambda_1), Z(\lambda_2)) = \mbox{cov} (Z(g \lambda_1), Z(g \lambda_2)), \quad \lambda_1, \lambda_2 \in \Lambda_{\kappa}, g \in \mathcal{G}.
\]
Therefore, $Z(x)$ is an IRF$\kappa$.

``$\Rightarrow$", assume that $Z(x)$ is an IRF$\kappa$, $\kappa \ge 1$. For any $y \in S^2$, consider the following measure
\[
\lambda_y(dx) = \delta_y(dx) - \sum_{l < \kappa} \sum_{m=-l}^l Y_l^m(y) Y_l^m(dx),
\]
where $\delta_y(dx)$ is the Dirac measure. First, for $l_0 < \kappa, |m_0| \le l_0$,
\begin{eqnarray*}
&\,& Y_{l_0}^{m_0}(\lambda_y) = \int_{S^2} Y_{l_0}^{m_0}(x) \left( \delta_y(dx) - \sum_{l < \kappa} \sum_{m=-l}^l Y_l^m(y) Y_l^m(dx) \right) \\
&\,& \quad = Y_{l_0}^{m_0}(y) - \sum_{l<\kappa} \sum_{m=-l}^l Y_l^m(y) \int_{S^2} Y_{l_0}^{m_0}(x) Y_l^m(dx).
\end{eqnarray*}
By orthogonality of the associated Legendre polynomials, the second term in the last equation reduces to $Y_{l_0}^{m_0}(y)$. Therefore, $Y_{l_0}^{m_0}(\lambda_y)=0$, that is, $\lambda_y$ annihilates all associated Legendre polynomials with order less than $\kappa$, which implies
\[
\lambda_y(dx) \in \Lambda_{\kappa}.
\]
For $Z(x)$, we have 
\[
Z(\lambda_y) = Z(y) - \sum_{l<\kappa} \sum_{m=-l}^l Z_{l,m} Y_l^m(y),
\]
which is exactly the low-frequency truncated process at $y$, i.e., 
\[
Z(\lambda_y)  = Z_{\kappa}(y).
\]
In the same way, 
\[
Z(g\lambda_y) = Z_{\kappa} (gy), g \in \mathcal{G}.
\] 
Therefore, for any $x, y \in S^2, g \in \mathcal{G}$,  we have $\mbox{cov} (Z_{\kappa}(x), Z_{\kappa}(y)) =   \mbox{cov} (Z(\lambda_x), Z(\lambda_y))$ and  $\mbox{cov} (Z_{\kappa}(gx), Z_{\kappa}(gy)) = \mbox{cov} (Z(g \lambda_x), Z(g \lambda_y))$. Note that  $Z(x)$ is an IRF$\kappa$ with $\mbox{cov} (Z(\lambda_x), Z(\lambda_y)) = \mbox{cov} (Z(g \lambda_x), Z(g \lambda_y))$ for $\lambda_x, \lambda_y \in \Lambda_{\kappa}$, we obtain
\[
\mbox{cov} (Z_{\kappa}(x), Z_{\kappa}(y)) = \mbox{cov} (Z_{\kappa}(gx), Z_{\kappa}(gy)).
\]
Similarly, for $\lambda \in \Lambda_{\kappa}$,
\[
\mbox{E} (Z_{\kappa}(x)) = \mbox{E} (Z_{\kappa}(gx)).
\]
That is, $Z_{\kappa}(x)$ is stationary.
\end{proof}

When $\kappa=0$, the IRF$0$ is the same as a stationary process. This is slightly different from \cite{Matheron1973}, where IRF$(-1)$ is a stationary process. For notation simplification, we assume $\kappa \ge 1$ throughout the rest of this manuscript. Huang et al. \cite{Huangetal2016} introduced IRFs on the circle and showed that an IRF on the circle is characterized by its low-frequency truncated process. Theorem 1 shows that the parallel result holds on the sphere and we further conjecture that this is true for hyper-spheres. \\

\noindent
{\bf Remark 1.} In Euclidean spaces, IRFs are associated with differential operations \cite{Matheron1973}. For example, the $(\kappa+1)$th derivative of a differentiable IRF$\kappa$ on the real line is stationary (\cite{ChilesandDelfiner2012}, Chapter 4). Theorem 1 indicates that low-frequency truncation operation of an IRF$\kappa$ on the sphere becomes stationary.  The implication of this result for splines on the sphere is discussed in Section 2.4. \\

Let $Z(x)$ be an IRF$\kappa$, then its low-frequency truncated process $Z_{\kappa}(x)$ is stationary (Theorem 1).  We have \cite{Obukhov1947, Yaglom1961}
\[
\mbox{cov} (Z_{l,m}, Z_{l',m'}) = a_l I(l, l') I(m,m'), \quad l,l' \ge \kappa, |m| \le l, |m'| \le l'.
\]
Hence, the covariance of $Z_{\kappa}(\cdot)$ is
\[
\mbox{cov} (Z_{\kappa}(x), Z_{\kappa}(y)) = \sum_{l=\kappa}^\infty a_l \sum_{m=-l}^l Y_l^m(x) Y_l^m(y). 
\]
By the addition theorem of Legendre polynomials, we arrive at
\begin{equation} \label{eq:al}
\mbox{cov} (Z_{\kappa}(x), Z_{\kappa}(y)) = \sum_{l=\kappa}^\infty \frac{2l+1}{4\pi} a_l P_l(\cos (d(x,y))).
\end{equation}
This covariance function plays an essential role in universal kriging on the sphere in Section 2.2. In addition, it offers the base in constructing reproducing kernel Hilbert space (Section 2.3). We denote this function by $\phi(\cdot)$, and name it the {\it intrinsic covariance function} of an IRF$\kappa$ on the sphere.\\

\noindent
{\bf Remark 2.} A more abstract view of IRFs in Euclidean spaces is based on equivalence class (\cite{ChilesandDelfiner2012}, Chapter 4).  Here, we extend this discussion to IRFs on the sphere.  Let $Z(x)$ be an IRF$\kappa$ and $A_{l,m}, l < \kappa, |m| \le l$ be random variables, a new random process
\begin{equation} \label{eq:eq1}
Z^*(x) = Z_{\kappa}(x) + \sum_{l < \kappa} \sum_{m=-l}^l A_{l,m} Y_l^m(x),
\end{equation}
is also an IRF$\kappa$, where $Z_{\kappa}(x)$ is the low-frequency truncated process of $Z(x)$. In fact, for $\lambda \in \Lambda_{\kappa}$, we have $Y_l^m(\lambda) = 0$ for $l < \kappa, |m| \le l$ and
\[
Z^*(\lambda) = Z_{\kappa}(\lambda) + \sum_{l < \kappa} \sum_{m=-l}^l A_{l,m} Y_l^m(\lambda) = Z_{\kappa}(\lambda).
\]
Therefore, $Z^*(x)$ and $Z(x)$ are both IRF$\kappa$ and they share the same truncated process $Z_{\kappa}(x)$. All random functions sharing the same truncated process $Z_{\kappa}(x)$ form an equivalence class. A natural choice of the representation of this equivalence class is $Z_{\kappa}(x)$, and all others are of form of (\ref{eq:eq1}). \\

\subsection{Universal kriging}

One of the primary applications of developing IRFs is universal kriging. Assume that $Z(x)$ is an IRF$\kappa$, its intrinsic covariance function is $\phi(\cdot)$, and the mean function is 
\begin{equation} \label{eq:betalm}
\mbox{E} (Z(x)) = \sum_{l < \kappa} \sum_{m=-l}^l \beta_{l,m} Y_l^m(x), 
\end{equation}
where $\beta_{l,m}$ are the coefficients.  Suppose that the data
\begin{equation} \label{eq:model}
\{ (x_i, w_i), i=1,\ldots,n \}, \quad n > (2\kappa-1),  \quad x_1, \ldots, x_n \in S^2
\end{equation}
are observed from this process with uncorrelated measurement errors
\[
W(x) = Z(x) + \epsilon(x), \quad x \in S^2,
\]
where $\epsilon(\cdot)$ is white noise with mean zero and variance $\sigma^2$.

To obtain the best linear unbiased estimator at $x_0 \in S^2$, following \cite{Cressie1993, Huangetal2016}, we start with a linear estimator
\[
\hat{Z}(x_0) = \utilde{\eta}^T \utilde{w}, 
\]
where $\utilde{w}=(w_1,\ldots, w_n)^T$ is the data vector and $\utilde{\eta}=(\eta_1, \ldots, \eta_n)^T$ is to be determined. For notation simplicity, we rearrange the lower spherical harmonics $\{Y_l^m(x), l < \kappa, |m| \le l\}$  and denote them by $\{q_{\nu}(x), \nu=1,\ldots, 2\kappa-1\}$. The unbiasedness leads to
\[
\utilde{\eta}^T Q = \utilde{q}^T,
\]
where
\[
Q=\{ q_{\nu}(x_i) \}_{n \times (2\kappa-1)}, \quad \utilde{q}=(q_1(x_0), \ldots, q_{2\kappa-1}(x_0))^T.
\]
Hence,
\[
\sum_{i=1}^n \eta_i q_j(x_i) = q_j(x_0), \quad j=1,\ldots, 2\kappa-1,
\]
which implies that the discrete measure
\[
\left( \sum_{i=1}^n \eta_i \delta(x_i) \right) - \delta(x_0) \in \Lambda_{\kappa}.
\]
The squared prediction error can be shown to be
\[
\mbox{E} \{ \hat{Z}(x_0) - Z(x_0) \}^2 = \sigma^2 \utilde{\eta}^T \utilde{\eta} + \utilde{\eta}^T \Psi \utilde{\eta} - 2 \utilde{\eta}^T \utilde{\phi} + \phi(0),
\]
where
\[
\Psi = \{ \phi(d(x_i, x_j)) \}_{n \times n}, \quad \utilde{\phi} = (\phi(d(x_1,x_0)),\ldots, \phi(d(x_n,x_0)))^T.
\]
The goal of universal kriging is to minimize the squared prediction error, subject to the unbiasedness constraints. Let $\utilde{\rho}$ be a Lagrange-multiplier vector of size $(2\kappa-1)$, we minimize the following with respect to both $\utilde{\rho}$ and $\utilde{\eta}$
\[
M(\utilde{\eta}, \utilde{\rho}) = \sigma^2 \utilde{\eta}^T \utilde{\eta} + \utilde{\eta}^T \Psi \utilde{\eta} - 2 \utilde{\eta}^T \utilde{\phi} + \phi(0) + 2 (\utilde{\eta}^T Q-\utilde{q}^T) \utilde{\rho}.
\]
One, then, can derive the universal kriging formula as
\begin{equation} \label{eq:uk}
\left\{ \begin{array}{l} (\Psi + \sigma^2 I) \utilde{\eta} + Q \utilde{\rho} = \utilde{\phi}, \\ Q^T \utilde{\eta} = \utilde{q}. \end{array} \right.
\end{equation}
It is clear from this kriging formula that there is no need to estimate the coefficients $\{ \beta_{l,m}\}$ in (\ref{eq:betalm}), showing the advantage of universal kriging \cite{Matheron1973}. Note that when $\kappa=1$,  universal kriging reduces to ordinary kriging \cite{Cressie1993} and the intrinsic covariance function $\phi(\cdot)$ relates to the semi-variogram directly, see \cite{Huangetal2011, Huangetal2016}. 

\subsection{Splines and kriging}

The goal of kriging is to find a linear unbiased predictor at an un-sampled location. It can also be viewed as surface prediction (\cite{Cressie1993}, Chapter 3), which has close ties to smoothing splines \cite{Wahba1990a}. The connection between splines and kriging have been extensively discussed in the literature, including \cite{Matheron1981, Watson1984, Lorenc1986, Cressie1989, Cressie1990, Cressie1993, Wahba1990a, Wahba1990b, KentandMardia1994, Laslett1994, FurrerandNychka2007, Huangetal2016}. Furrer and Nychka \cite{FurrerandNychka2007} showed that, given a covariance function, one can construct a reproducking kernel and obtain a general spline estimate in Euclidean spaces. 

Based on \cite{Levesleyetal1999}, Huang et al. \cite{Huangetal2016} established the connection between IRF and RKHS, and provided the corresponding smoothing formula for universal kriging on the circle. Now, we extend this discussion to the sphere. Note that, Levesley et al. \cite{Levesleyetal1999} constructed the RKHS based on a non-negative integer $\kappa^*$ and a sequence $\{ \beta_l^* \}_{l=\kappa}^\infty$. While both $\kappa^*$ and $\beta_l^*$ play essential roles in \cite{Levesleyetal1999}, their motivation and interpretation are unclear. In this section, we demonstrate that $\kappa^*$ is exactly the order of an IRF, and $\beta_l^*$ are closely related to the IRF's intrinsic covariance fucntion $\phi(\cdot)$.  

Now, following \cite{Levesleyetal1999}, we formally introduce RKHS. Based on the intrinsic covariance function $\phi(\theta)$ of an IRF$\kappa$, if $\phi(\theta)$ is of the form (\ref{eq:al}), we define a function space
\[
X_{\kappa} = \left\{ f(x) = \sum_{l=0}^\infty \sum_{m=-l}^l a_{l,m} Y_l^m(x), x \in S^2: \quad \sum_{l=\kappa}^\infty \frac{1}{a_l} \sum_{m=-l}^l |a_{l,m}|^2 < \infty \right\}.
\]
Then, for $f, g \in X_{\kappa}$, a semi-inner product is defined as
\begin{equation} \label{eq:semi}
\langle f, g \rangle_{\kappa} = \sum_{l=\kappa}^\infty \frac{1}{a_l} \sum_{m=-l}^l a_{l,m,f} a_{l,m,g}.
\end{equation}
There is a nil space for this semi-inner product
\begin{equation} \label{eq:nil}
N = \mbox{span} \{ Y_l^m(x), 0 \le l < \kappa, m=-l,\ldots,l \},
\end{equation}
with $\mbox{dim} (N) = d_N = 2\kappa-1$. Let $\{ \tau_1, \ldots, \tau_{d_N} \} \in S^2$ be a set of distinct points. Then, the inner product 
\[
\langle f, g \rangle = \sum_{\nu=1}^{d_N} f(\tau_{\nu}) g(\tau_{\nu}) + \langle f, g \rangle_{\kappa}
\]
is well defined and $X_{\kappa}$ can be shown to be complete \cite{Levesleyetal1999} with respect to the norm induced by the inner product (\ref{eq:semi}). In addition, there exist $p_1(\cdot), \ldots, p_{d_N}(\cdot) \in N$ such that $p_{\nu}(\tau_{\mu}) = I(\nu,\mu)$ for $1 \le \nu,\mu \le d_N$. As discussed in \cite{Levesleyetal1999}, the space $X_{\kappa}$ is a Hilbert space in which point evaluations are continuous linear functionals. Therefore, for $x,y \in S^2$, there exists a reproducing kernel 
\begin{eqnarray*}
&\,& H(x,y) = \phi(d(x,y)) - \sum_{\nu=1}^{d_N} \{ \phi(d(x,\tau_{\nu})) p_{\nu}(y) + \phi(d(y,\tau_{\nu})) p_{\nu}(x) \} \\
&\,& \quad + \sum_{\nu=1}^{d_N} \sum_{\mu=1}^{d_N} \phi(d(\tau_{\nu}, \tau_{\mu})) p_{\nu}(x) p_{\mu}(y) + \sum_{\nu=1}^{d_N} p_{\nu}(x) p_{\nu}(y).
\end{eqnarray*}

Hence, given the intrinsic covariance function of an IRF$\kappa$, we construct a RKHS. Conversely, it is obvious that $H(x,y)$, being a reproducing kernel, is positive definite. Given this property and from the Kolmogorov existence theorem, there exists a Gaussian random process with $H(x,y)$ as its covariance function. Suppose we construct a random process $Z(x)$ whose covariance function is $H(x,y)$. We now show that $Z(x)$ is an IRF$\kappa$ with the intrinsic covariance function $\phi(d(x,y))$ for $x, y \in S^2$. Note that measures in $\Lambda_{\kappa}$ annihilate any functions in nil space $N$, that is, for $\lambda \in \Lambda_{\kappa}$, $p_{\nu}(\lambda) = 0, \nu=1, \ldots, d_N$. Therefore, 
\[
\mbox{cov} (Z(\lambda_1), Z(\lambda_2)) = \int_{S^2} \int_{S^2} \phi(d(x,y)) \lambda_1(dx) \lambda_2(dx), \quad \lambda_1, \lambda_2 \in S^2,
\]
This covariance is rotation invariant, hence, $Z(x)$ is an IRF$\kappa$ with $\phi(\cdot)$ as its intrinsic covariance function.

We can see from this development that $\kappa^*$ in \cite{Levesleyetal1999} is exactly the order of an IRF, and $\beta^{*2}_l = 1/a_l$. Such interpolation helps understand the RKHS theory on the sphere in \cite{Levesleyetal1999}. Furthermore, we will show that the smoothing formula in \cite{Levesleyetal1999} finds its counterpart in universal kriging. Given observed data $\{ (x_i, w_i) \}_{i=1}^n$ in (\ref{eq:model}), the smoothing problem is to find a function $f_{\alpha}(x) \in X_{\kappa}$ such that it minimizes \cite{Taijeronetal1994, Levesleyetal1999}
\begin{equation} \label{eq:sf}
\sum_{i=1}^n (w_i - f(x_i))^2 + \alpha \| f \|^2_{\kappa},
\end{equation}
where $\alpha > 0$ is the smoothing parameter and $\| \cdot \|_{\kappa}$ is induced by the semi-inner  product (\ref{eq:semi}). The minimizer can be shown to be
\begin{equation} \label{eq:smoothingformula}
f_{\alpha}(x) = \sum_{\nu=1}^{d_N} b_{\nu} q_{\nu}(x) + \sum_{i=1}^n c_i \phi(d(x_i,x)),
\end{equation}
where $\utilde{b} = (b_1, \ldots, b_{d_N})^T$ and $\utilde{c}=(c_1, \ldots, c_n)^T$ satisfying
\begin{equation} \label{eq:eq2}
\left\{ \begin{array}{l} (\Psi+\alpha I)\utilde{c} + Q \utilde{b} = \utilde{w}, \\ Q^T \utilde{c} = \utilde{0}, \end{array} \right.
\end{equation}
where $\Psi$ and $Q$ are the same as given in (\ref{eq:uk}). 
To show the connection between this smoothing formula and universal kriging, note that the smoothing formula (\ref{eq:smoothingformula}) at an unobserved point $x_0 \in S^2$ can be re-written in the following manner
\begin{eqnarray*}
&\,& f_{\alpha}(x_0) = (\utilde{c}^T,\utilde{b}^T) \left( \begin{array}{c} \utilde{\phi} \\ \utilde{q} \end{array} \right),
\end{eqnarray*}
where $\utilde{\phi}$ and $\utilde{q}$ are the same in (\ref{eq:uk}). Solve for $(\utilde{c}^T, \utilde{b}^T)$ in (\ref{eq:eq2}) and we have
\[
f_{\alpha}(x_0) = (\utilde{w}^T, \utilde{0}^T_{d_N \times1}) \left( \begin{array}{cc} \Psi + \alpha I & Q \\ Q^T & 0 \end{array} \right)^{-1} \left( \begin{array}{c} \utilde{\phi} \\ \utilde{q} \end{array} \right).
\]
Therefore, $f_{\alpha}(x_0)$ can be re-written as
\[
f_{\alpha}(x_0)=(\utilde{w}^T, \utilde{0}^T) \left( \begin{array}{c} {\utilde{\eta}}^* \\ {\utilde{\rho}}^* \end{array} \right),
\]
where
\[
\left( \begin{array}{c} {\utilde{\eta}}^* \\ {\utilde{\rho}}^* \end{array} \right) = \left( \begin{array}{cc} \Psi + \alpha I & Q \\ Q^T & 0 \end{array} \right)^{-1} \left( \begin{array}{c} \utilde{\phi} \\ \utilde{q} \end{array} \right),
\]
or
\begin{equation} \label{eq:smoothing}
\left\{ \begin{array}{l} (\Psi+\alpha I) {\utilde{\eta}}^* + Q {\utilde{\rho}}^* = \utilde{\phi}, \\ Q^T {\utilde{\eta}}^* = \utilde{q}. \end{array} \right.
\end{equation}

\noindent
{\bf Remark 3.} These two equations (\ref{eq:uk}) and (\ref{eq:smoothing}) are exactly the dual formulas of universal kriging \cite{Cressie1993, ChilesandDelfiner2012}, but expressed using the intrinsic covariance function of an IRF on the sphere introduced in this paper. Usually, universal kriging is viewed as a linear estimator of observed data and the smoothing formula as the sum of lower spherical harmonic trends and the linear combination of the intrinsic covariance. Based on equations (\ref{eq:uk}) and (\ref{eq:smoothing}), these two views are essentially the same. The connection between universal kriging and the smoothing formula on the sphere is, thus, obvious. In so doing, we find that the lower monomials in Euclidean spaces need to be replaced by lower spherical harmonics. This echoes the low-frequency truncation that replaced the differential operations in Euclidean spaces, see {Remark 1}.\\

\noindent
{\bf Remark 4.} The smoothing parameter $\alpha$ in (\ref{eq:smoothing}) equals the noise variance $\sigma^2$ in (\ref{eq:uk}), and plays the same role in prediction.  For example, in the smoothing formula (\ref{eq:smoothingformula}), when $\alpha$ increases to infinity, the minimization procedure demands $\| f \|_{\kappa}$ approaches zero, which shows that $\utilde{c} \to 0$, and the smoothing formula reduces to a spherical harmonic regression. In kriging, when $\sigma^2 \to \infty$, the noise overwhelms the spatial dependency, the process becomes essentially uncorrelated. The squared prediction error is dominated by $\sigma^2 \utilde{\eta}^T \utilde{\eta}$. Then, universal kriging reduces to minimize $\utilde{\eta}^T \utilde{\eta}$ subject to the unbiasedness restriction, which also leads to the same spherical harmonic regression prediction. When both $\alpha$ and $\sigma^2$ decrease to zero, both smoothing and kriging result in exact interpolation. 

\subsection{Splines on the sphere and thin-plate splines}

Splines on the sphere are commonly used for surface fitting \cite{Wahba1981, Wahba1990a, Robeson1997}. Given the data $\{(x_i, w_i), i=1,\ldots, n\}$ in (\ref{eq:model}),  Wahba \cite{Wahba1981} introduced splines on the sphere as a function to minimize the following
\[
\sum_{i=1}^n (w_i - f(x_i))^2 + \alpha J(f),
\]
where the penalty term is through the Laplacian operator $\Delta$ on the sphere and for $m$ even,
\[
J_m(f) = \int_{S^2} \left( \Delta^{m/2} f \right)^2 \sin \zeta d \zeta d \psi.
\]
The spline estimator is shown to be \cite{Wahba1981}
\begin{equation} \label{eq:wahba81spline}
f_{\alpha}(x) = d + \sum_{i=1}^n c_i K(x,x_i),
\end{equation}
where
\[
K(x,y) = \frac{1}{4\pi} \sum_{l=1}^\infty \frac{2l+1}{l^m(l+1)} P_l (\cos (d(x,y))), \quad x, y \in S^2.
\]
In Euclidean spaces, the integer $m$ plays a significant role and relates to the order of IRFs in kriging \cite{KentandMardia1994}. While $K(x,y)$ here takes the place of $\phi(\cdot)$ in (\ref{eq:smoothingformula}), the integer $m$ only alters the form of $K(x,y)$ and loses its connection to the order $\kappa$ in IRF on the sphere. Furthermore, while the nil space in (\ref{eq:sf}) is spanned by the lower spherical harmonics of orders up to $\kappa$, splines on the sphere in (\ref{eq:wahba81spline}) only include a constant term. By viewing splines on the sphere through the smoothing formula (\ref{eq:smoothingformula}), we now see that splines on the sphere act simply as ordinary kriging.  This observation echoes Remark 1 and reveals that splines with derivative penalty have limited applications on the sphere (Huang et al. \cite{Huangetal2016} observed the similar result on the circle). Then, a more appropriate smoothing (spline) model on the sphere is the smoothing formula (\ref{eq:sf}) in Section 2.3.  A general approach for hyper-spheres can be found in \cite{Taijeronetal1994}.

For surface fitting in multi-dimensional Euclidean spaces, thin-plate splines are often used \cite{Duchon1977,  Meiguet1979, Wahba1990a}. Some applications \cite{Newetal2002} have extended thin-plate splines to  the sphere by mapping longitude and latitude directly on two-dimensional Euclidean spaces. It shall be warned that such practice is flawed. The radial basis kernel  of thin-plate splines are  \cite{Duchon1977, Meiguet1979}
\[
E(x,y) = E(d^*(x,y)) = |d^*(x,y)|^2 \log |d^*(x,y)|,
\]
where $d^*(x,y)$ is the Euclidean distance between $x$ and $y$. When thin-plate splines are obtained in Euclidean spaces, $E(\cdot)$ is shown \cite{Wahba1990a} to be conditionally positive definite of order $2$. However, when $x, y$ are on the sphere, the distance between $x$ and $y$ is the spherical distance $d(x,y)$ in (\ref{eq:dxy}). To obtain a valid thin-plate spline on the sphere, $E(d(x,y))$ needs to be conditionally positive definite of order $2$. Menengatto and Peron \cite{MegengattoandPeron2004} showed that the coefficients of $E(d(x,y))$ in the Legendre polynomials expansion must be non-negative for all orders higher than or equal to $2$. However, direct computation shows that $E(\cdot)$ yields negative coefficients for such expansion on the sphere. For example, the fourth-order coefficient is negative. Hence, thin-plate splines are not directly applicable for surface fitting on the sphere. 

\subsection{Extensions}


An IRF$\kappa$ is characterized by its truncated process, while the mean function falls in nil space (\ref{eq:nil}). This is a direct extension of IRFs in Euclidean spaces. It is important and noteworthy that the spherical harmonics are orthogonal, while the monomials used in Euclidean spaces are not necessarily so. Therefore, in developing the notion of IRFs, the allowable measure can also be defined to annihilate the spherical harmonics of orders in a set $\{i_1, \ldots, i_{\kappa} \}$ of integers, instead of  just $\{1, \ldots, \kappa \}$. Therefore, paralleling to Theorem 1, a continuous random process in the sense of quadratic mean is an IRF$(i_1,\ldots,i_{\kappa})$ on the sphere if and only if
\[
Z_{(i_1,\ldots,i_{\kappa})}(x) = \sum_{l \notin \{i_1,\ldots,i_{\kappa}\}} \sum_{m=-l}^l Z_{l,m} Y_l^m(x)
\]
is stationary. This is a more general result that includes the IRF$\kappa$ developed in Section 2.1 as a special case and may offer more flexibility in modeling the mean function of an IRF on the sphere.



\bibliography{mybibfile}{}

\begin{thebibliography}{10}

\bibitem{Aronszajn1950}
N~Aronszajn.
\newblock Theory of reproducing kernels.
\newblock {\em Transactions of the American Mathematical Society},
  68:337--€"404, 1950.

\bibitem{ChilesandDelfiner2012}
J.~Chil\`es and P.~Delfiner.
\newblock {\em Geostatistics: Modeling Spatial Uncertainty}.
\newblock Wiley, New York, 2 edition, 2012.

\bibitem{Cressie1989}
N.~Cressie.
\newblock Geostatistics.
\newblock {\em The American Statistician}, 43:197--€"202, 1989.

\bibitem{Cressie1990}
N.~Cressie.
\newblock Reply to ``comment on cressie'' by {G}. {W}ahba.
\newblock {\em The American Statistician}, 44:256--€"258, 1990.

\bibitem{Cressie1993}
N.~Cressie.
\newblock {\em Statistics for Spatial Data, revised ed.}
\newblock Wiley, New York, 1993.

\bibitem{Duchon1977}
J.~Duchon.
\newblock Splines minimizing rotation-invariant semi-norms in {S}obolev spaces.
\newblock In {\em Constructive Theory of Functions of Several Variables}, pages
  85--€"100. Springer-Verlag, Berlin, 1977.

\bibitem{FurrerandNychka2007}
E.~Furrer and D.~Nychka.
\newblock A framework to understand the asymptotic properties of kriging and
  splines.
\newblock {\em Journal of the Korean Statistical Society}, 36:57--€"76, 2007.

\bibitem{GelfandandVilenkin1964}
I.M. Gel'fand and N.Y. Vilenkin.
\newblock {\em Applications of Harmonic Analysis (Generalized Functions)},
  volume~4.
\newblock Academica Press, New York, 1964.

\bibitem{Gneiting2013}
T.~Gneiting.
\newblock Strictly and non-strictly positive definite functions on spheres.
\newblock {\em Bernoulli}, 19:1327--€"1349, 2013.

\bibitem{Heatonetal2014}
M.~Heaton, M.~Katzfuss, C.~Berrett, and D.~Nychka.
\newblock Constructing valid spatial processes on the sphere using kernel
  convolutions.
\newblock {\em Environmetrics}, 25:2--€"15, 2014.

\bibitem{Huangetal2011}
C.~Huang, H.~Zhang, and S.~Robeson.
\newblock On the validity of commonly used covariance and variogram functions
  on the sphere.
\newblock {\em Mathematical Geosciences}, 43:721--€"733, 2011.

\bibitem{Huangetal2012}
C.~Huang, H.~Zhang, and S.~Robeson.
\newblock A simplified representation of the covariance structure of axially
  symmetric processes on the sphere.
\newblock {\em Statistics and Probability Letters}, 82:1346--1351, 2012.

\bibitem{Huangetal2016}
C.~Huang, H.~Zhang, and S.~Robeson.
\newblock Intrinsic random functions and universal kriging on the circle.
\newblock {\em Statistics and Probability Letters}, 108:33--€"39, 2016.

\bibitem{Levesleyetal1999}
Levesley J., W.~Light, D.~Ragozin, and X.~Sun.
\newblock A simple approach to the variational theory for interpolation on
  spheres.
\newblock {\em International Series of Numerical Mathematics}, 132:117--€"143,
  1999.

\bibitem{Jones1963}
A.~H. Jones.
\newblock Stochastic processes on a sphere.
\newblock {\em Annals of Mathematical Statistics}, 34:213--€"217, 1963.

\bibitem{KentandMardia1994}
J.~Kent and K.~Mardia.
\newblock The link between kriging and thin-plate splines.
\newblock In F.P. Kelly, editor, {\em Probability, Statistics, and
  Optimization: A tribute to Peter Whittle}, pages 325--€"339. Wiley,
  Chichester, 1994.

\bibitem{Laslett1994}
G.F. Laslett.
\newblock Kriging and splines: An empirical comparison of their predictive
  performance in some applications.
\newblock {\em Journal of American Statistics Association}, 89:391--€"400,
  1994.

\bibitem{Lorenc1986}
A.C. Lorenc.
\newblock Analysis methods for numerical weather prediction.
\newblock {\em Quarterly Journal of the Royal Meteorological Society},
  112:1177--€"1194, 1986.

\bibitem{MarinucciandPeccati2011}
D.~Marinucci and G.~Peccati.
\newblock {\em Random Fields on the Sphere: Representation, Limit Theorems, and
  Cosmological Applications}.
\newblock London Mathematical Society, Lecture Notes Series 389, Cambridge,
  2011.

\bibitem{Matheron1973}
G.~Matheron.
\newblock The intrinsic random functions and their applications.
\newblock {\em Advances in Applied Probability}, 5:439--€"468, 1973.

\bibitem{Matheron1979}
G.~Matheron.
\newblock Comment translater les catastrophes. {L}a structure des {F.A.I}.
  \'{g}en\'{e}rales. manuscript.
\newblock Technical Report N-167, Centre de G\'{e}ostatistique, Fontainebleau,
  France, 1979.

\bibitem{Matheron1981}
G.~Matheron.
\newblock Splines and kriging: their formal equivalence.
\newblock In {\em Syracuse University Geological Contributions}. Syracuse
  University, 1981.

\bibitem{Meiguet1979}
J.~Meiguet.
\newblock Multivariate interpolation in arbitrary points made simple.
\newblock {\em Journal of Applied Mathematics Physics}, 30:292--€"304, 1979.

\bibitem{MegengattoandPeron2004}
V.~Menegatto and A.~Peron.
\newblock Conditionally positive definite kernels on {E}uclidean domains.
\newblock {\em Journal of Mathematical Analysis and Applications},
  294:345--€"359, 2004.

\bibitem{Newetal2002}
M.~New, D.~Lister, M.~Hulme, and I.~Makin.
\newblock A high-resolution data set of surface climate over global land areas.
\newblock {\em Climate Research}, 21:1--€"25, 2002.

\bibitem{Obukhov1947}
A.~M. Obukhov.
\newblock Statistical homogeneous random fields on a sphere.
\newblock {\em Uspekhi Matematicheskikh}, 2:196--€"198, 1947.

\bibitem{Robeson1997}
S.M. Robeson.
\newblock Spherical methods for spatial interpolation: Review and evaluation.
\newblock {\em Cartography and Geographic Information Systems}, 24:3--€"20,
  1997.

\bibitem{Roy1969}
R.~Roy.
\newblock {\em Processus stochastiques sur la sph\`ere}.
\newblock PhD thesis, Universit\'e de Montr\'eal, 1969.

\bibitem{Roy1973}
R.~Roy.
\newblock Estimation of the covariance function of a homogeneous process on the
  sphere.
\newblock {\em Annals of Statistics}, 1:780--€"785, 1973.

\bibitem{Royden1988}
H.L. Royden.
\newblock {\em Real Analysis}.
\newblock Prentice Hall, Englewood Cliffs, NJ, 3rd edition, 1988.

\bibitem{Schoenberg1942}
I.J. Schoenberg.
\newblock Positive definite functions on spheres.
\newblock {\em Duke Mathematics Journal}, 9:96--€"108, 1942.

\bibitem{Stein2007}
M.~L. Stein.
\newblock Spatial variation of total column ozone on a global scale.
\newblock {\em Annals of Applied Statistics}, 1:191--€"210, 2007.

\bibitem{Stein1999}
M.L. Stein.
\newblock {\em Statistical Interpolation of Spatial Data: Some Theory for
  Kriging}.
\newblock Springer, New York, 1999.

\bibitem{Taijeronetal1994}
H.J. Taijeron, A.G. Gibson, and C.~Chandler.
\newblock Spline interpolation and smoothing on hyperspheres.
\newblock {\em SIAM Journal on Scientific Computing}, 15:1111--€"1125, 1994.

\bibitem{Wahba1981}
G.~Wahba.
\newblock Spline interpolation and smoothing on the sphere.
\newblock {\em SIAM Journal on Scientific and Statistical Computing}, 2, 1981.

\bibitem{Wahba1990b}
G.~Wahba.
\newblock Comment on {C}ressie.
\newblock {\em The American Statistician}, 44:255--€"256, 1990.

\bibitem{Wahba1990a}
G.~Wahba.
\newblock {\em Spline Models for Observational Data}.
\newblock CBMS-NSF regional conference series in applied mathematics,
  Philadelphia, 1990.

\bibitem{Watson1984}
G.S. Watson.
\newblock Smoothing and interpolation by kriging with splines.
\newblock {\em Mathematical Geology}, 16:601--€"615, 1984.

\bibitem{Yaglom1955}
A.M. Yaglom.
\newblock Correlation theory of processes with random stationary nth increments
  (in {R}ussian).
\newblock {\em Matematicheskii Sbornik}, 37:141--€"196, 1955.
\newblock English translation in (1958) {\it American Mathematical Society
  Translations: Series 2}, {\bf 8}, 87-141. American Mathematical Society,
  Providence, R. I., 1958.

\bibitem{Yaglom1961}
A.M. Yaglom.
\newblock Second-order homogeneous random fields.
\newblock {\em Fourth Berkeley Symposium on Mathematical Statistics and
  Probability}, 2:593--€"622, 1961.

\bibitem{ZhangandHuang2014}
H.~Zhang and C.~Huang.
\newblock A note on processes with random stationary increments.
\newblock {\em Statistics and Probability Letters}, 94:153--€"161, 2014.

\bibitem{ZhuandWu2010}
Z.~Zhu and Y.~Wu.
\newblock An efficient algorithm for estimation and prediction of a class of
  convolution-based spatial nonstationary models.
\newblock {\em Journal of Computational and Graphical Statistics}, 19:74--€"95,
  2010.

\end{thebibliography}

\end{document}